\documentclass[reqno, 12pt,a4paper]{amsart}
\usepackage[utf8]{inputenc}	
\usepackage[T1]{fontenc}
\usepackage[in,headings]{fullpage}
\usepackage{amsbsy, amscd, amsfonts, amsmath, amsrefs, amssymb, amsthm}
\usepackage{graphicx}
\usepackage{float}
\usepackage{color}   
\usepackage[colorlinks=true, pdfstartview=FitV, linkcolor=blue, 
            citecolor=blue, urlcolor=blue, hyperfigures=true, bookmarksnumbered]{hyperref}
\usepackage{comment}
\excludecomment{A}
\usepackage{algorithm} 
\usepackage{algpseudocode} 

\newtheorem{theorem}{Theorem}
\newtheorem*{cor}{Corollary}

\newtheorem{proposition}[theorem]{Proposition}
\newtheorem{lemma}[theorem]{Lemma}

\theoremstyle{definition}
\newtheorem{definition}[theorem]{Definition}
\newtheorem{remark}[theorem]{Remark}

\theoremstyle{remark}

\newcommand{\C}{\mathbf{C}}
\newcommand{\Z}{\mathbf{Z}}

\renewcommand{\Re}{\mathop{\mathrm{Re}}\nolimits}
\renewcommand{\Im}{\mathop{\mathrm{Im}}\nolimits}
\newcommand{\Rzeta}{\mathop{\mathcal R }\nolimits}

\newfont{\cmbsy}{cmbsy10}
\newfont{\cmmib}{cmmib10}
\newcommand{\Orden}{\mathop{\hbox{\cmbsy O}}\nolimits}

\begin{document}

\title[Zeros on the fourth quadrant]
{Zeros of $\Rzeta(s)$ on the fourth quadrant.}
\author[Arias de Reyna]{J. Arias de Reyna}
\address{%
Universidad de Sevilla \\ 
Facultad de Matem\'aticas \\ 
c/Tarfia, sn \\ 
41012-Sevilla \\ 
Spain.} 




\email{arias@us.es, ariasdereyna1947@gmail.com}


\begin{abstract}
We show that there is a sequence of zeros of $\Rzeta(s)$ in the fourth quadrant. We show that the $n$-th zero $\rho_{-n}=\beta_{-n}+i\gamma_{-n}$, with $\beta_{-n}\sim 4\pi^2 n/\log^2n$ and $\gamma_{-n}\sim-4\pi n/\log n$. We give the first terms of an asymptotic development of $\rho_{-n}$ and an algorithm to calculate $\rho_{-n}$ from $n$.
\end{abstract}

\maketitle
\section{Introduction}
Numerical evidence suggests that the auxiliary function of Riemann $\Rzeta(s)$ \cite{A166}, has three main lines of zeros.
\begin{figure}[H]
\begin{center}
\includegraphics[width=0.8\hsize]{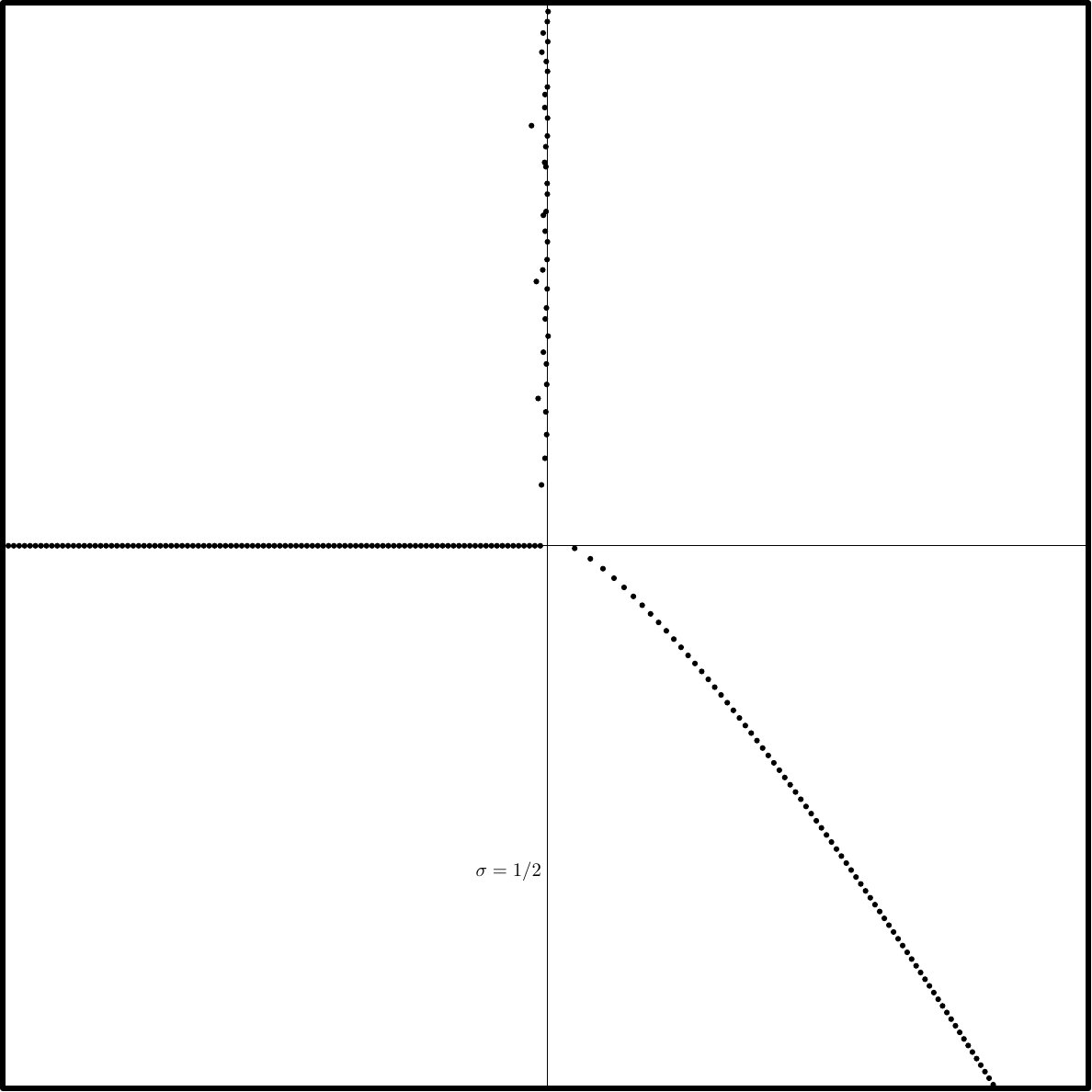}
\caption{Zeros of $\Rzeta(s)$ on the rectangle $(-200,200)^2$.}
\label{default}
\end{center}
\end{figure}
One line along the negative real axis is formed with the trivial zeros at $s=-2n$ for $n=1$, $2$, \dots. The more interesting zeros follow approximately the positive imaginary axis; some of his properties are studied in \cite{A173}, \cite{A98}, and numerically in 
\cite{A172}.  Here we try to get some information about the zeros on the fourth quadrant.

In Section \ref{S:2} we give a parametric form \eqref{linezeros} of the curve defined by the zeros. This is achieved because this line of zeros is contained in the limit of two asymptotic expansions given in \cite{A100} that are valid at each side of the line.

In Section \ref{S:3} we obtain an approximate equation for the zeros in the third quadrant. First, the main equation $\zeta(s)=\Rzeta(s)+\chi(s)\overline{\Rzeta}(1-s)$ implies that  $\rho\notin\Z$ is a zero of $\Rzeta(s)$ if and only if the symmetric $\tilde \rho$ with respect to the critical line satisfies the equation $\chi(1-\tilde \rho)\Rzeta(\tilde \rho)=\zeta(1-\tilde\rho)$. The zero is at a position where no asymptotic is known, but in $\tilde\rho$ there is a nice approximation. This is what allows us to obtain approximate equations. The equation for a zero $f(\eta)=1$ is approximated by a simple function $g(\eta)$, with $f(\eta)=g(\eta)(1+\Orden(|\eta|^{-1})$ (see Proposition \ref{facts}). Here $\eta$ is an auxiliar variable connected to $s$.

In Section \ref{S:4} we give a general theorem  showing that two functions related by $f(\eta)=g(\eta)(1+\Orden(|\eta|^{-1})$ and satisfying some mild conditions will have very close solutions of equations $f(\eta)=1$ and $g(\eta)=1$. In Section \ref{S:5} we apply this theorem to relate the zeros of $\Rzeta(s)$ in the fourth quadrant to the zeros of $g(\eta)=1$. We have to introduce a simpler approximate equation $h(\eta)=1$ that can be solved using the Lambert $W$ function.

Section \ref{S:6} contains an algorithm to calculate the $n$-th zero $\rho_{-n}$ for $n\ge0$ of $\Rzeta(s)$ in the fourth quadrant starting from the value $n$.  We have used this algorithm to obtain the first 2122 zeros with 25 correct decimal digits. The table of zeros is included after the \texttt{end\{document\}} line of the \TeX\  file associated with this document.

In Section \ref{S:7} we give the first terms of the asymptotic expansion  of $\rho_{-n}$.

\section{Location of Zeros} \label{S:2}

\begin{proposition}\label{P:1}
There is some $R>0$ such that for any zero $\rho$ of $\Rzeta(s)$ in the fourth quadrant with $|\rho|>R$,  there are some $r>0$ and $0<\epsilon<1$ such that 
\begin{equation}\label{linezeros}
\rho=2\pi r^2\exp\Bigl\{i\Bigl(-\frac{\pi}{2}+\frac{\pi}{2\log r}-\frac{\pi^3}{24\log^3 r}+\frac{\pi^3\epsilon}{24\log^4r}\Bigr)\Bigr\}.
\end{equation}
\end{proposition}
\begin{proof}
Let $\rho $ be a zero of $\Rzeta(s)$ in the fourth quadrant. Let $r=\sqrt{|\rho|/2\pi}$ and $\rho=2\pi r^2 e^{i\theta}$ with $-\frac{\pi}{2}<\theta<0$. By \cite{A100}*{Thm.~16} and \cite{A100}*{Prop.~6} the point $\rho\notin M$. It follows (see \cite{A100}*{Remark  17}) that
for $|\rho|>2\pi e^2$
\[-\frac{\pi}{2}+\frac{\pi}{2\log r}-\frac{\pi^3}{24\log^3r}+\frac{\pi^5}{160\log^5r}+\Orden(\log^{-6}r)<\theta.\]
Therefore, if $|\rho|>R$ we will have $\Re\rho>2$. 
By \cite{A100}*{Cor.~14} $\rho\notin L$. Therefore, if $|\rho|>R$ we have (see Remark 12 and Theorem 12 in \cite{A100})
\[\theta<-\frac{\pi}{2}+\frac{\pi}{2\log r}-\frac{\pi^3}{24\log^3 r}+\frac{\pi^3}{48\log^4r}+\Orden(\log^{-5}r).\]
So for $|\rho|>R$ with $R$ large enough, we will have 
\[-\frac{\pi}{2}+\frac{\pi}{2\log r}-\frac{\pi^3}{24\log^3 r}<\theta<-\frac{\pi}{2}+\frac{\pi}{2\log r}-\frac{\pi^3}{24\log^3 r}+\frac{\pi^3}{24\log^4r}.\]
It follows that there exists $0<\epsilon<1$ satisfying \eqref{linezeros}.
\end{proof}

\section{Approximate equation for zeros}\label{S:3}

Recall that we denote by  $\chi(s):=2^s\pi^{s-1}\sin\frac{\pi s}{2}\Gamma(1-s)$  the meromorphic function appearing in the functional equation $\zeta(s)=\chi(s)\zeta(1-s)$.

\begin{proposition}\label{condition}
For any number $s\in\C$ let $\tilde s$ be the symmetric respect the critical line $\sigma=\frac12$, that is $\tilde s:=1-\overline{s}$. Let $\rho\in\C\smallsetminus \Z$, the two following equations are equivalent
\[\Rzeta(\rho)=0,\quad \chi(1-\tilde \rho)\Rzeta(\tilde \rho)=\zeta(1-\tilde\rho).\]
\end{proposition}
\begin{proof}
Siegel \cite{Siegel}*{Eq.~(56)} (but see also \cite{A166}*{Eq.~(10)})  proved that 
\[\zeta(s)=\Rzeta(s)+\chi(s)\overline{\Rzeta}(1-s),\qquad \text{with } \overline{\Rzeta}(s)=\overline{\Rzeta(\overline{s})}. \]
If we assume that $\Rzeta(\rho)=0$, then putting $s=\tilde\rho$ in this equation, we get $\overline{\Rzeta}(1-\tilde\rho)=\overline{\Rzeta(\rho)}=0$. Therefore, $\zeta(\tilde\rho)=\Rzeta(\tilde\rho)$. Multiplying by $\chi(1-\tilde\rho)$ we get $\chi(1-\tilde \rho)\Rzeta(\tilde \rho)=\zeta(1-\tilde\rho)$.

We may reverse the reasoning except when $\chi(1-\tilde\rho)=0$ or $\chi(1-\tilde\rho)=\infty$. But the zeros of $\chi(s)$ are $s=0$, $-2$, $-4$, \dots\, and its poles are $s=1$, $s=3$, \dots are all integers.
\end{proof}

If $\rho$ is a zero of $\Rzeta(s)$ in the fourth quadrant, we do not have an asymptotic expansion valid for $\Rzeta(s)$ at this point. But we have one valid for $\Rzeta(s)$ for $s=\tilde\rho$. 

\begin{proposition}\label{P:defregion}
We define  the following  region of the complex plane
\begin{equation}\label{defOmega}
\Omega=\Bigl\{re^{i\phi}\colon e<r, \frac{\pi}{4}<\phi<\frac{\pi}{2}-\frac{\log r}{r^2}\Bigr\}.
\end{equation}
For $\eta\in\Omega$ we have $\zeta(-2\pi i \eta^2)=1+\Orden(|\eta|^{-1})$.
\end{proposition}
\begin{proof}
For $\eta\in\Omega$ we have $-2\pi i \eta^2=2\pi r^2 e^{-\frac{\pi i}{2}+2i\phi}$ with the argument $0<-\frac{\pi}{2}+2\phi<\frac{\pi}{2}-\frac{2\log r}{r^2}$. So, it is in the first quadrant. Its real part is 
\[\Re(-2\pi i \eta^2)=2\pi r^2\cos(-\tfrac\pi2+2\phi)\ge 2\pi r^2 \sin(\pi-\tfrac{2\log r}{r^2})>12,\quad\text{ for $r>e$}.\]  It follows that $\zeta(-2\pi i\eta^2)=1+\Orden(2^{2\pi i \eta^2})$.
And we have
\begin{align*}
|2^{2\pi i\eta^2}|&=\exp(\Re(2\pi i\eta^2)\log2)\le\exp\bigl(-2\pi r^2\log2\sin(\tfrac{2\log r}{r^2})\bigr)\\ &\le \exp(-(\pi\log2)\log r)<\frac{1}{r}.\qedhere  \end{align*}
\end{proof}

\begin{proposition}\label{P:mainaprox}
For $\eta\in\Omega$ we have 
\begin{equation}\label{firstaprox}
\Rzeta(1+2\pi i\eta^2)=-\chi(1+2\pi i\eta^2)\eta^{2\pi i\eta^2}e^{-\pi i \eta^2}\frac{\sqrt{2} e^{\frac{3\pi i}{8}}\sin(\pi\eta)}{2\cos(2\pi\eta)}(1+\Orden(|\eta|^{-1})).
\end{equation}
\end{proposition}
\begin{proof}
For $\eta\in\Omega$ we have $\eta=re^{i\phi}$ with $r>e$ and $\frac{\pi}{4}<\phi<\frac{\pi}{2}$. Therefore, $2\pi i \eta^2$ is in the third quadrant and $s=1+2\pi i \eta^2$ is in the set $N$ considered in \cite{A100}*{Thm.~29}. Therefore, with the notations in that paper (\cite{A100}*{eq.~(79) and (80)})
\begin{equation}\label{E:left}
\Rzeta(s)=\frac{(-1)^m}{2i}
\chi(s)\eta^{2\pi i\eta^2}e^{-\pi i \eta^2}\overline{G}(p)(1+\Orden(|\eta|^{-1})).
\end{equation}
By definition here $\eta=\sqrt{\frac{s-1}{2\pi i}}$, taking the root so that its argument is in the interval $(-\pi/4,3\pi/4)$, therefore, coincides with our initial $\eta$. The numbers
$\eta_1$ and $\eta_2$ are the real and imaginary parts of $\eta$.
Then $m=\lfloor\eta_1+\eta_2\rfloor$ and 
\[p=-2(m+\tfrac12-(\eta_1+i\eta_2))=2\bigl\{(\eta_1+\eta_2-\lfloor\eta_1+\eta_2\rfloor-\tfrac12)-(1-i)\eta_2\bigr\}.\]
So, with $\lambda$ a real number with $-1\le\lambda\le 1$ we have 
\[p=\lambda-2\sqrt{2}e^{-\pi i/4}\eta_2.\]
The function $G(q)$ is defined in \cite{A100}*{Eq.~(13)}
\[G(q):=\frac{e^{\frac{\pi i}{2}q^2}-\sqrt{2}e^{\frac{\pi i}{8}}\cos\frac{\pi q}{2}}{\cos\pi q},\quad \overline{G}(p)=\overline{G(\overline p)}=\frac{e^{-\frac{\pi i}{2}p^2}-\sqrt{2}e^{-\frac{\pi i}{8}}\cos\frac{\pi p}{2}}{\cos\pi p}.\]
\[\overline{G}(p)=-\frac{\sqrt{2}e^{-\frac{\pi i}{8}}\cos\frac{\pi p}{2}}{\cos\pi p}\Bigl(1-\frac{e^{-\frac{\pi i}{2}p^2}}{\sqrt{2}e^{-\frac{\pi i}{8}}\cos\frac{\pi p}{2}}\Bigr).\]
\[-\frac{\pi i}{2}p^2=-\frac{\pi i}{2}(\lambda-2\sqrt{2}e^{-\pi i/4}\eta_2)^2=
-\frac{\pi i}{2}(\lambda^2-4\sqrt{2}\lambda e^{-\pi i/4}\eta_2-8i\eta_2^2).\]
\[|e^{-\frac{\pi i}{2}p^2}|= e^{2\pi\lambda\eta_2-4\pi\eta_2^2}\le e^{-4\pi\eta_2^2+2\pi\eta_2}.\]
\[|\cos\tfrac{\pi p}{2}|\ge\sinh(|\Im(\pi p)/2|)=\sinh(\pi\eta_2).\]
Since $\eta=r e^{i\phi}$ with $\frac{\pi}{4}<\phi<\frac{\pi}{2}$ it follows that $\eta_2>\frac{1}{\sqrt{2}}|\eta|>e/\sqrt{2}$, it follows that 
\[\Bigl|\frac{e^{-\frac{\pi i}{2}p^2}}{\sqrt{2}e^{-\frac{\pi i}{8}}\cos\frac{\pi p}{2}}\Bigr|\le \frac{e^{-4\pi\eta_2^2+2\pi\eta_2}}{\sqrt{2}\sinh(\pi \eta_2)}\ll e^{-2\pi\eta_2^2}.\]
Therefore, 
\[\Rzeta(1+2\pi i \eta^2)=-\frac{(-1)^m}{2i}\chi(1+2\pi i\eta^2)\eta^{2\pi i\eta^2}e^{-\pi i \eta^2}\frac{\sqrt{2}e^{-\frac{\pi i}{8}}\cos\frac{\pi p}{2}}{\cos\pi p}(1+\Orden(|\eta|^{-1})).\]
Since $m$ is an integer, we have
\[\cos\pi p=\cos(-2\pi m-\pi+2\pi\eta)=-\cos(2\pi\eta),\quad \cos\frac{\pi p}{2}=\cos(\pi m+\tfrac{\pi}{2}-\pi\eta)=(-1)^m\sin\pi\eta.\]
This proves \eqref{firstaprox}.
\end{proof}

\begin{definition}\label{deffunctions}
We will consider several functions defined for $\Im(\eta)>0$ with values in $\C$
\begin{align*}
f(\eta)&:=\frac{\chi(-2\pi i\eta^2 )\Rzeta(1+2\pi i\eta^2)}{\zeta(-2\pi i \eta^2)},\\
u(\eta)&:=2\eta^2\log\eta-\eta^2+\eta+i\frac{\log2}{2\pi}-\frac{1}{8},&
v(\eta)&:=2\eta^2\log\eta-\eta^2+i\frac{\log2}{2\pi}-\frac{1}{8},\\
g(\eta)&:=e^{\pi i u(\eta)},&
h(\eta)&:=e^{\pi i v(\eta)},
\end{align*}
where $\log\eta$  denotes the main branch of the logarithm with $0<\Im\log\eta<\pi$. 
\end{definition}

\begin{proposition}
Let $\eta\in\Omega$ be such that $f(\eta)=1$, then $\rho=2\pi i\overline\eta^2$ is a zero of $\Rzeta(s)$ in the fourth quadrant.
\end{proposition}
\begin{proof}
Let $\rho=2\pi i \overline\eta^2$. Since $\eta=re^{i\phi}\in\Omega$ we have 
$\rho=2\pi r^2 e^{\frac{\pi i}{2}-2i\phi}$, with $0>\frac{\pi}{2}-2\phi>-\frac{\pi}{2}+\frac{2\log r}{r^2}$ and $r>e$. Therefore, $\rho$ is in the fourth quadrant. We have \[\overline\rho=-2\pi i\eta^2,\quad \tilde\rho=1+2\pi i\eta^2,\quad 1-\tilde\rho=-2\pi i\eta^2.\]
Then $f(\eta)=1$ implies \[\chi(1-\tilde\rho)\Rzeta(\tilde\rho)=\zeta(1-\tilde\rho).\] By Proposition \ref{condition}, this is equivalent to  $\Rzeta(\rho)=0$.
\end{proof}

\begin{proposition}\label{Rzetaf}
There is some $R\ge 2\pi e^2$ such that if $\rho$ be a zero of $\Rzeta(s)$ in the fourth quadrant with $|\rho|>R$, then $\rho=2\pi i \overline\eta^2$, where $\eta\in\Omega$ is such that $f(\eta)=1$.
\end{proposition}
\begin{proof}
Let $\rho$ be a zero of $\Rzeta(s)$ in the fourth quadrant. By Proposition \ref{P:1} it can be put in the form \eqref{linezeros}. Then define 
\begin{equation}\label{rhoform}
\eta=\sqrt{\frac{i\overline\rho}{2\pi}}=r\exp\Bigl\{i\Bigl(\frac{\pi}{2}-\frac{\pi}{4\log r}+\frac{\pi^3}{48\log^3 r}-\frac{\pi^3\epsilon}{48\log^4 r}\bigr)\Bigr\}.
\end{equation}
Since $R>2\pi e^2$, we have $|\eta|>e$.
It is clear that there is some constant $R$ such that if $|\rho|=2\pi r^2>R$ then this $\eta\in\Omega$. And by Proposition \ref{condition} we will have $f(\eta)=1$. 

A numerical study shows that $R=2\pi e^2$ is a good election, but we will not show this. In any case, there are 8 zeros of $\Rzeta(s)$ in the fourth quadrant with $|\rho|<2\pi e^2$. 
\end{proof}

\begin{proposition}\label{facts}
For $\eta\in\Omega$, we have  $f(\eta)=g(\eta)(1+\Orden(|\eta|^{-1}))$.
\end{proposition}
\begin{proof}
In Proposition \ref{P:mainaprox} we have proved that for $\eta\in \Omega$ we have \eqref{firstaprox} and by Proposition \ref{P:defregion} $\zeta(-2\pi i\eta^2)=1+\Orden(|\eta|^{-1})$. Therefore, 
\[f(\eta)=\chi(-2\pi i\eta^2)\cdot\Bigl(-\chi(1+2\pi i\eta^2)\eta^{2\pi i\eta^2}e^{-\pi i \eta^2}\frac{\sqrt{2} e^{\frac{3\pi i}{8}}\sin(\pi\eta)}{2\cos(2\pi\eta)}(1+\Orden(|\eta|^{-1}))\Bigr).\]
The function $\chi(s)$ satisfies $\chi(s)\chi(1-s)=1$, so that 
\begin{equation}\label{Siegeleq}
f(\eta)=-\eta^{2\pi i\eta^2}e^{-\pi i \eta^2}\frac{\sqrt{2} e^{\frac{3\pi i}{8}}\sin(\pi\eta)}{2\cos(2\pi\eta)}(1+\Orden(|\eta|^{-1}))\Bigr).
\end{equation}
As mentioned above for $\eta\in\Omega$ we have $\eta_2>\frac{1}{\sqrt{2}}|\eta|$. Therefore,
\[\sin(\pi\eta)=-\frac{1}{2i}e^{-\pi i \eta}(1-e^{2\pi i\eta})=\frac{i}{2}e^{-\pi i \eta}(1+\Orden(e^{-2\pi \eta_2}))=\frac{i}{2}e^{-\pi i \eta}(1+\Orden(|\eta|^{-1})).\]
In the same way 
\[2\cos(2\pi\eta)=e^{-2\pi i \eta}(1+\Orden(|\eta|^{-1})).\]
It follows that 
\begin{align*}
f(\eta)&=-\eta^{2\pi i\eta^2}e^{-\pi i \eta^2}\frac{\sqrt{2} e^{\frac{3\pi i}{8}}\frac{i}{2}e^{-\pi i \eta}}{e^{-2\pi i \eta}}(1+\Orden(|\eta|^{-1}))\\
& =\exp\bigl(2\pi i\eta^2\log\eta-\pi i\eta^2+\pi i \eta-\tfrac{\log 2}{2}-\tfrac{\pi i}{8}\bigr)(1+\Orden(|\eta|^{-1}))\\
&=g(\eta)(1+\Orden(|\eta|^{-1})).\qedhere
\end{align*}
\end{proof}

\section{X-rays of similar functions}\label{S:4}
We first prove a general theorem. This theorem explains why the X-rays of two functions with the same asymptotic behavior almost coincide (under certain conditions).  Applying it to the functions $f$ and $g$ of Definition \ref{deffunctions}, we reduce the study of the complicated equation $f(\eta)=1$ to the simpler equation $g(\eta)=1$.
\begin{theorem}\label{X-rayTheorem}
Let $V\subset \C$ be an open set of the complex plane. Assume that for each $r>0$ the open set $\Omega_r:=\{z\in V\colon |z|>r\}$ is  non empty. 
Let $f$, $g\colon V\to\C$ be holomorphic functions such that:
\begin{itemize}
\item[(a)] $f(z)=g(z)(1+\Orden(|z|^{-1}))$.
\item[(b)] There are constants $0<c<C$, sucht that $c|z|\log|z|\le |g'(z)/g(z)|\le C |z|\log|z|$, for  all $z\in V$.
\item[(c)] Let $A=\{a\in V\colon f(a)=1\}$ and $B=\{b\in V\colon g(b)=1\}$. There exist two positive constants $r_0$ and $\rho_0$ such that if $w\in A\cup B$ satisfies $|w|>r_0$, then $\overline{D}(w, \rho_0)\subset V$.
\end{itemize}
Then there is $R$ such that for any $b\in B$ with $|b|>R$ there is a unique $a\in A$ such that $|b-a|<|b|^{-2}$. The induced map $b\mapsto a$ is injective, and its image contains all $a\in A$ with $|a|>R+1$.
\end{theorem}

\begin{proof}
Let $b\in B$, with $|b|>R>1$, by hypothesis (c), taking $R>r_0$ the disc $\overline{D}(b; \rho_0)\subset \Omega$. 
To prove the main assertion we will apply Rouche's Theorem to the functions $f(z)-1$ and $g(z)-1$ on the circle $|z-b|=|b|^{-2}$. Taking $R$ large enough so that $R^{-1}<\rho_0$. First, we need a bound of $g(z)$ on a larger disc.  
For $|z-b|\le \rho_0$ we have 
\[\log|g(z)|=\log|g(z)|-\log|g(b)|=\Re \int_b^z\frac{g'(u)}{g(u)}\,du,\]
integrating along the segment $[b,z]$. Therefore, for $|z-b|\le |b|^{-3/2}$ 
\[\log|g(z)|\le \int_b^z\Bigl|\frac{g'(u)}{g(u)}\Bigr|\,|du|\le C(|b|+|b|^{-3/2})\log(|b|+|b|^{-3/2})\cdot|z-b|\le C'|b|\log|b|\cdot|z-b|.\]
On $|z-b|=|b|^{-3/2}$ we have 
\[|g(z)|\le\exp(C'|b|\log|b|\cdot|z-b|)\le \exp(C'|b|^{-1/2}\log|b|)<2,\]
taking $R$ big enough. With this bound for $g(z)$,  Cauchy's inequality for the coefficients of a power series yields 
\[g(z)=g(b)+g'(b)(z-b)+\sum_{n=2}^\infty a_n(z-b)^n, \qquad |a_n|\le 2|b|^{3n/2}.\]
Therefore, for $|z-b|\le |b|^{-2}$ and $R$ large enough
\begin{align*}
|g(z)-1|&=|g(z)-g(b)|\ge|g'(b)(z-b)|-2|z-b|\sum_{n=2}^\infty|b|^{3n/2} |z-b|^{n-1}\\
&\ge \Bigl(c|b|\log|b| -2|b|^{2}\frac{1}{|b|(1-|b|^{-1/2})}\Bigr)|z-b|\ge \bigl(c|b|\log|b|-4|b|\bigr)|z-b|\\ &> (c'|b|\log|b|) |z-b|.
\end{align*}
It follows that $z=b$ is the unique zero of $g(z)-1$ in the disc $|z-b|\le |b|^{-3/2}$, and is a simple zero.

By the above reasoning  for  $|z-b|= |b|^{-2}$ we have 
\begin{equation}\label{boundinf}
|g(z)-1|>
\frac{c'\log|b|}{|b|}.
\end{equation}  
We are in position to apply Rouche's Theorem to $f(z)-1$ and $g(z)-1$ on the disc $|z-b|\le |b|^{-2}$. We have for $|z-b|=|b|^{-2}$  
\begin{align*}
|(f(z)-1)-(g(z)-1)|&=|f(z)-g(z)|\le \frac{B}{|z|}|g(z)|, && \text{by hypothesis (a)}\\
&\le \frac{2B}{|b|-|b|^{-2}}\le \frac{4B}{|b|} \\
&\le \frac{c'\log|b|}{|b|},  && \text{taking $R$ large enough}\\
& < |g(z)-1|, && \text{by \eqref{boundinf}.}
\end{align*}
Therefore, Rouche's Theorem implies that $f(z)-1$ and $g(z)-1$ have the same number of zeros in the disc $|z-b|\le|b|^{-2}$. By \eqref{boundinf} the function $g(z)-1$ vanishes only at $b$, and this is simple zero because $g'(b)\ne0$. 

We have seen that for any point  $b\in B$ with $|b|>R$,  there is a unique point $a\in A$ with $|b-a|\le |b|^{-2}$.

To see that the correspondence is injective and covers the large zeros of $f(z)-1$, we notice that the above reasoning is only local near a point $b\in B$. The roles of $f$ and $g$ can be reversed locally, as we see below.  

Condition (a) implies that $g(z)=f(z)(1+\Orden(|z|^{-1}))$, since $(1+\Orden(|z|^{-1}))^{-1}$ is equal to $(1+\Orden(|z|^{-1}))$.

Condition (a) implies $\frac{f'(z)}{f(z)}=\frac{g'(z)}{g(z)}+\frac{u'(z)}{1+u(z)}$ where 
$u(z)\le B/|z|$ for $z\in \Omega_r$  for some fixed $r$, (take $R$ greater than this $r$). If $f(a)=1$ and $|a|>R+1$, by hypothesis (c), the disc with center $a$ and radius $\rho_0$ is contained in $\Omega_r$ and then for $|z-a|<\rho_0/2$, we have (assuming $|a|>R>2\rho_0$)
\[|u'(z)|=\Bigl|\frac{1}{2\pi i}\int_{|\zeta-a|=\rho_0}\frac{u(\zeta)}{(\zeta-z)^2}\,d\zeta\Bigr|\le \frac{B}{|a|-\rho_0}\frac{2\pi \rho_0}{2\pi\rho_0^2}\le \frac{C}{|a|}\]
It follows that taking $|a|>R$ with $R$ large enough
\[c'|z|\log |z|\le c|z|\log|z|-\frac{C}{|a|}\frac{1}{1-B/|z|}\le \Bigl|\frac{f'(z)}{f(z)}\Bigr|\le \Bigl|\frac{g'(z)}{g(z)}\Bigr|+\frac{C}{|a|}\frac{1}{1-B/|z|}\le C'|z|\log|z|.\]
This is only locally around each point $a$, but this is sufficient to apply Rouche's Theorem to the two discs $|z-a|\le 4|a|^{-2}$ and $|z-a|\le \frac14|a|^{-2}$.  Therefore, for $|a|>R+1$ we get a unique $b$ with $g(b)=1$ and $|a-b|\le 4|a|^{-2}$ and then this same $b$ satisfies $|a-b|\le\frac14|a|^{-2}$. 

Then $|b|\ge |a|-|b-a|\ge R+1-\frac14|a|^{-2}>R$. Therefore, there is a unique $a'\in A$ with 
$|b-a'|\le |b|^{-2}$. But we know that $|a-b|\le\frac14|a|^{-2}$. Therefore, $a=a'$ if we have 
$\frac14|a|^{-2}<|b|^{-2}$. This is equivalent to $|b|\le2|a|$. This is true since 
$|b|\le |a|+|b-a|\le |a|+\frac14|a|^{-2}\le 2|a|$, because $|a|>R>1$.
\end{proof}

\begin{cor}
In the hypothesis of the Proposition, all the zeros of $f(z)-1=0$ in $\Omega_R$ are simple.
\end{cor}
Because condition (c) implies that $g'$ does not vanish on the set $\Omega$. A similar reasoning applies to $f$ locally near each zero.

\section{The approximate equation}\label{S:5}
\begin{figure}[H]
\begin{center}
\includegraphics[width=0.9\hsize]{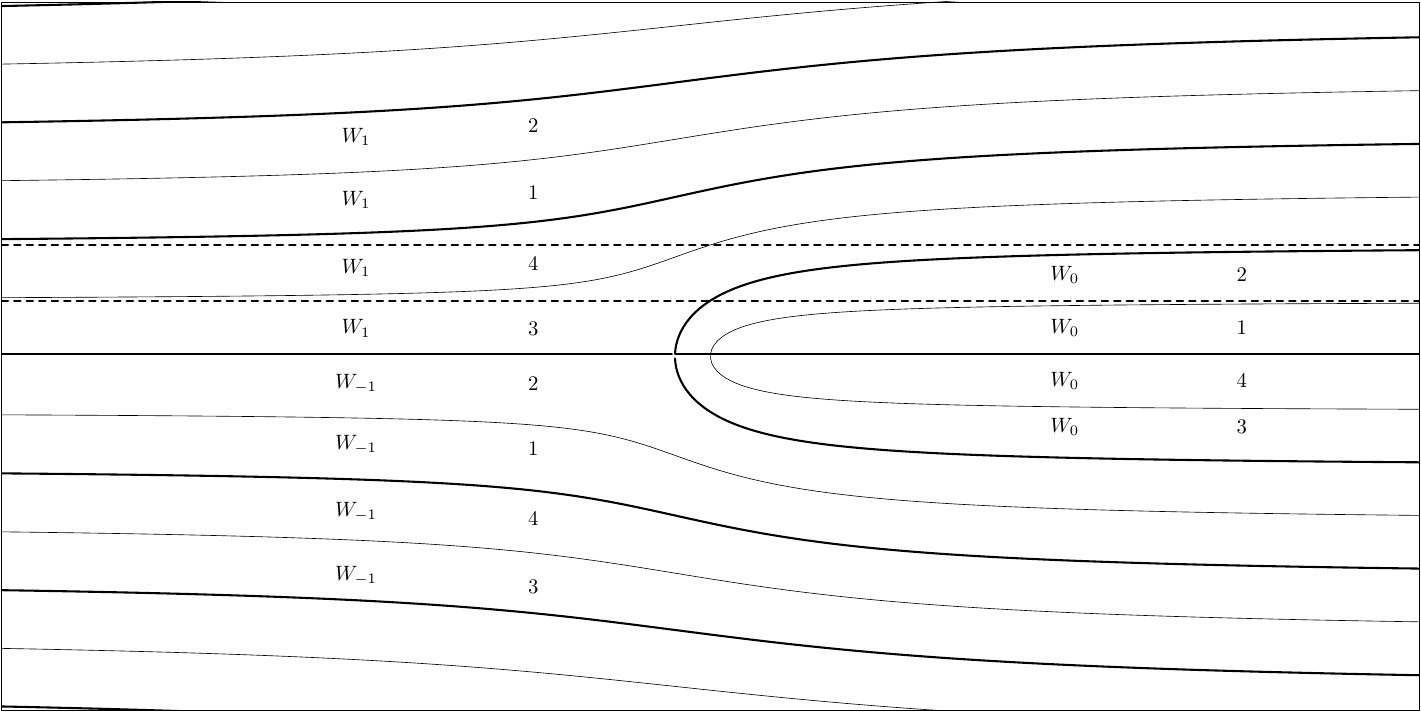}
\caption{X-ray of the function $z e^z$ in $(-20,20)\times(-10,10)$.}
\label{X-ray}
\end{center}
\end{figure}

Before using Theorem \ref{X-rayTheorem} we need some information about the points $\eta\in\Omega$ with $g(\eta)=1$. First, we consider the points where $h(\eta)=1$. 
\begin{proposition}\label{givesn}
For each point $\eta\in\Omega$  with $h(\eta)=1$ there is a natural number $n$ such that 
\begin{equation}\label{eqWeta}
\log\frac{\eta^2}{e}=W_1\Bigl(\frac{1}{8e}-\frac{2n}{e}-i\frac{\log2}{2\pi e}\Bigr).
\end{equation}
where $W_1$ is the first branch of the Lambert W function (see \cite{C}).
\end{proposition}
\begin{proof}
By definition $h(\eta)=e^{\pi i v(\eta)}$. Therefore, if $h(\eta)=1$ there is an integer $n$ such that $v(\eta)=-2n$. This equation is equivalent to 
\begin{equation}\label{E:toasymp}
\frac{\eta^2}{e}\log\frac{\eta^2}{e}=\frac{1}{e}\Bigl(\frac18-2n-i\frac{\log2}{2\pi}\Bigr).
\end{equation}
The Lambert $W$ function is the inverse of $z\mapsto z e^z$. Therefore, the above equation is equivalent $\log(\eta^2/e)=W(\frac{1}{8e}-\frac{2n}{e}-i\frac{\log2}{2\pi e})$. The function $W$ as the logarithm has many branches. They are defined in \cite{C}. The X-ray of the function $ze^z$ (that is, the lines where this function is real or purely imaginary; see Figure \ref{X-ray}) divides the plane into regions that are transformed in a quadrant. The derivative of $ze^z$ only vanishes at $z=-1$ where the function is real. It follows that each of the regions defined by the X-ray transform bijectively in a quadrant (see \cite{C}*{p.~344}). In Figure \ref{X-ray}, I have put a numeral $1$, $2$, $3$ or $4$ in each region, indicating the quadrant in which it transforms. We also put in each region a symbol $W_{-1}$, $W_0$ and $W_1$ indicating that the corresponding region is in the image by the corresponding branch of $W$. The real and imaginary  lines of the X-ray are asymptotic to the lines $\Im z=k\frac{\pi}{2}$.  

If $\eta=re^{i\phi}\in\Omega$, then $\pi/4<\phi<\pi/2$ and $r>e$. Then $\log(\eta^2/e)=2\log r-1+2i\phi$ is contained between the two dashed lines $\Im z=\pi/2$ and $\Im z=\pi$. Its image by $z e^z$ is $\frac18-2n-i\frac{\log2}{2\pi}$ is in the third of fourth quadrant. These two conditions implies that $\log(\eta^2/e)$ is the image of the branch $W_1$. Therefore, we have \eqref{eqWeta}.  

When $\eta\in \Omega$ the real part of $\eta^2/e$ is positive. Therefore, the image of $W_1$ is not only between the two dashed lines, but also has a positive real part. This is only possible if 
$\eta^2/e$ is the image of a point $e^{-1}(\frac18-2n-i\frac{\log2}{2\pi})$ in the third quadrant. This means $n\ge 1$.
\end{proof}

\begin{definition}
Let $\eta''_n$ be the number $\eta$ defined in \eqref{eqWeta}, taking the square root of $\eta^2$ with $-\pi/4\le \arg(\eta''_n)<3\pi/4$.
\end{definition}

\begin{proposition}\label{propeta''}
There is a natural number $N$ such that for $n\ge N$ we have $\eta''_n\in\Omega$. \begin{equation}\label{asymptotic}
\eta''_n=\eta_1+i\eta_2, \qquad \eta_1\sim\frac{\pi}{2\log n}\sqrt{\frac{2n}{\log n}},\quad
\eta_2\sim\sqrt{\frac{2n}{\log n}}.
\end{equation}
The circle with center at $\eta''_n$ and radius $2$ is contained in $\Omega$. 
\end{proposition}
\begin{proof}
Each branch of the Lambert $W$ function has an asymptotic expansion. We may obtain an asymptotic expansion of $\eta''_n$ when $n\to+\infty$ or $n\to-\infty$. We start with $n\to+\infty$.

In \cite{C}*{eq.~4.19} we find that
\[W(z)=L_1-L_2+\frac{L_2}{L_1}+\Orden(L_2^2/L_1^2),\]
where $L_1=\mathop{\rm Log} z$ and $L_2=\log\mathop{\rm Log} z$, $\log$ denoting the principal branch of the logarithm and $\mathop{\rm Log}$ any determination. In particular, the $W_k$ branch of the Lambert function is obtained taking $\mathop{\rm Log}z=2k\pi i+\log z$.

In our case when $n\to+\infty$ 
\begin{align*}
L_1&=2k\pi i+\log\Bigl(\frac{1}{8e}-\frac{2n}{e}-\frac{\log2}{2\pi e}i\Bigr)=2k\pi i+\log\Bigl\{-\frac{2n}{e}\Bigl(1-\frac{1}{16n}+\frac{\log2}{4\pi n}i\Bigr)\Bigr\}\\
&=\log n+(2k-1)\pi i +\log 2-1+\Orden(n^{-1}).
\end{align*}
Then 
\begin{align*}
L_2&=\log\bigl(\log n+(2k-1)\pi i +\log 2-1+\Orden(n^{-1})\bigr)\\
&=\log\log n+\log\Bigl(1+\frac{(2k-1)\pi i +\log 2-1}{\log n}+\Orden(1/n\log n)\Bigr)\\
&=\log\log n+\frac{(2k-1)\pi i +\log 2-1}{\log n}-\frac{((2k-1)\pi i +\log 2-1)^2}{2\log^2 n}+\Orden(\log^{-3}n).
\end{align*}
And then we get 
\[\frac{L_2}{L_1}=\frac{\log\log n}{\log n}\Bigl(1-\frac{(2k-1)\pi i +\log 2-1}{\log n}+\frac{(2k-1)\pi i +\log 2-1}{\log n\log\log n}+\Orden(\log^{-2}n)\Bigr).\]
Combining all this
\begin{align*}
W_k\Bigl(\frac{1}{8e}-\frac{2n}{e}&-i\frac{\log2}{2\pi e}\Bigr)=\log n-\log\log n+(2k-1)\pi i +\log 2-1+\frac{\log\log n}{\log n}\\&-\frac{(2k-1)\pi i +\log 2-1}{\log n}+\Orden((\log\log n/\log n)^2).
\end{align*}
This number will be $\eta^2/e$ with $\eta\in\Omega$ only if its imaginary part is in the interval $(\pi/2,\pi)$. For $n$ large, this is only true taking $k=1$.  

Now $\eta^2= e e^W$ and $\eta=e^{1/2} e^{W/2}$.
Therefore,
\begin{equation}\label{E:eta2}
{\eta''_n}^2=-\frac{2n}{\log n}\exp\Bigl(\frac{\log\log n}{\log n}-\frac{\log 2-1+\pi i}{\log n}\Bigr)\Bigl(1+\Orden((\log\log n/\log n)^2)\Bigr),
\end{equation}
and 
\[\eta''_n=i\sqrt{\frac{2n}{\log n}}\exp\Bigl(\frac{\log\log n}{2\log n}-\frac{\log 2-1+\pi i}{2\log n}\Bigr)\Bigl(1+\Orden((\log\log n/\log n)^2)\Bigr).\]
Expanding the exponential and retaining only terms smaller than  $(\log\log n/\log n)^{2}$ we get \begin{equation}\label{etan}
\eta''_n=i\sqrt{\frac{2n}{\log n}}\Bigl(1+\frac{\log\log n}{2\log n}-\frac{\log 2-1+\pi i}{2\log n}\Bigr)(1+\Orden((\log\log n/\log n)^2).
\end{equation}
This proves \eqref{asymptotic}. It follows that $\eta''_n=r e^{i\phi}$ with $\phi=\frac{\pi}{2}-\alpha$, where 
\[r\sim\sqrt{\frac{2n}{\log n}},\qquad \alpha=\frac{\pi}{2}-\arg(\eta''_n)\sim \arctan\frac{\pi}{2\log n}\sim \frac{\pi}{4\log r}.\]
The definition \eqref{defOmega} of $\Omega$ implies that there is some $N$ such that $\eta''_n\in\Omega$ for $n\ge N$.
A similar reasoning proves that $r e^{i\phi}(1+ \delta)\in\Omega$  for $|\delta|\le 2/r$. So that the disc $\overline D(\eta''_n,2)\subset \Omega$. 
\end{proof}

\begin{remark}
There is an asymptotic expansion for $\eta''_n$ in terms of $n$ of the form 
\begin{equation}\label{firstasymptotic}
\eta''_n=i\sqrt{\frac{2n}{\log n}}\Bigl(1+\sum_{k=1}^\infty\frac{U_k(\log\log n)}{\log^k n}\Bigr),
\end{equation}
where $U_k(x)$ is a polynomial of degree $k$.
\end{remark}

\begin{lemma}\label{L:InOmega}
Let $\eta\in\C$ with $\Im(\eta)>0$ and $|\eta|>e^2$. Then for $|\arg(\eta)-\frac{\pi}{4}|\le 1/10$ or $\frac{\pi}{2}-\frac{\log r}{r^2}\le\arg(\eta)\le\frac{\pi}{2}+\frac{1}{10}$, we have $\Im v(\eta)\ne 0$.
\end{lemma}
\begin{proof}
Let $\eta=re^{\frac{\pi i}{4}+i\phi}$, with $|\phi|\le \frac{1}{10}$. Then 
\begin{align*}
\Im(v(\eta))&=\Im\bigl(2r^2ie^{2i\phi}(\log r+\tfrac{\pi i}{4}+\phi i)-2r^2ie^{2i\phi}+i\tfrac{\log2}{2\pi}-\tfrac18\bigr)\\
&=2r^2(\log r\cos(2\phi)-(\tfrac{\pi}{4}+\phi)\sin(2\phi))-2r^2\cos(2\phi)+\tfrac{\log2}{\pi}\\
&=2r^2\log\tfrac{r}{e}\cos(2\phi)-2r^2(\tfrac{\pi}{4}+\phi)\sin(2\phi) +\tfrac{\log2}{\pi}\\
&\ge 2\cos(1/5)r^2\log \tfrac{r}{e}-2r^2(\tfrac{\pi}{4}+\tfrac{1}{10})\sin(1/5)+\tfrac{\log2}{\pi}>0,
\end{align*}
because $\cos(1/5)>(\pi/4+1/10)\sin(1/5)$.

If $\eta=re^{\frac{\pi i}{2}+i\phi}$, with $|\phi|\le \frac{1}{10}$, then 
\begin{align*}
\Im(v(\eta))&=\Im\bigl(-2r^2e^{2i\phi}(\log r+\tfrac{\pi i}{2}+\phi i)+2r^2e^{2i\phi}+i\tfrac{\log2}{2\pi}-\tfrac18\bigr)\\
&=-2r^2(\log r\sin(2\phi)+(\tfrac{\pi}{2}+\phi)\cos(2\phi))+2r^2\sin(2\phi)+\tfrac{\log2}{\pi}\\
&=-2r^2(\tfrac{\pi}{2}+\phi)\cos(2\phi)-2r^2\log \tfrac{r}{e}\sin(2\phi)+\tfrac{\log2}{\pi}.
\end{align*}
Therefore, since $\sin(2\varphi)>-2\log r/r^2$, we have
\begin{align*}
-\Im(v(\eta))&=2r^2(\tfrac{\pi}{2}+\phi)\cos(2\phi)+2r^2\log \frac{r}{e}\sin(2\phi)-\tfrac{\log2}{\pi}\\
&\ge\pi \cos(1/5) r^2-4\log r\log\tfrac{r}{e}-\tfrac{\log2}{\pi}>0.\qedhere
\end{align*}
\end{proof}

\begin{proposition}\label{nearaproxzeros}
There is a positive constant $R$ such that if $a\in\Omega$ with $g(a)=e^{\pi i u(a)}=1$, satisfies $|a|>R$, then there is some natural number $n$ such that $b=\eta''_n\in\Omega$ satisfies $b=a+\delta$ with $\delta=1/4\log a+\Orden(|a|^{-1})$.

There is no other point $a'\in\Omega$ with $u(a')=2n$.
\end{proposition}

\begin{proof}
Since $e^{\pi i u(a)}=0$, there is an integer $n$ with $u(a)=2n$. 
We apply Rouche's Theorem to $u(\eta)-2n$ and $v(\eta)-2n$ in the disc $|\eta-a|\le 1$, to show that there is a solution $b$ of $v(\eta)=2n$ with $|\delta|<1$. 
For  $|\eta-a|=1$, we have
\[|(u(\eta)-2n)-(v(\eta)-2n)|=|\eta|\le 1+|a|.\]
The derivatives are 
\[u'(\eta)=4\eta\log\eta+1,\quad u''(\eta)=4\log\eta+4,\quad u'''(\eta)=\frac{4}{\eta},\]
\[ u^{(k+2)}(\eta)=4(-1)^{k-1}\frac{(k-1)!}{\eta^k}.\]
Since $u(a)-2n=0$,  we have for $|\eta-a|=1$
\begin{align*}
&|u(\eta)-2n|=\\
&=\Bigl|(4a\log a+1)(\eta-a)+(2\log a+2)(\eta-a)^2+\sum_{k=1}^\infty\frac{4(-1)^{k-1}}{(k+2)(k+1)k a^k}(\eta-a)^{k+2}\Bigr|\\
&\ge 4|a||\log a|-1-2-2|\log a|-\frac{1}{|a|}>1+|a|,\quad \text{for $|a|>4$},
\end{align*}
since $\sum_{k=1}^\infty {4}/{(k+2)(k+1)k}=1$ and $|a|>1$. 

Therefore, for $|\eta-a|=1$ we have 
\[|(u(\eta)-2n)-(v(\eta)-2n)|\le 1+|a|<|u(\eta)-2n|.\]
And there is a zero $b$ of $v(\eta)-2n$ in the circle of center $a$ and radius $1$. 
We have $v(a+\delta)-2n=v(b)-2n=0$. Therefore, $u(a+\delta)-(a+\delta)-2n=0$. Expanding $u(a+\delta)-2n$ in Taylor series, we have 
\[-a-\delta+(4a\log a+1)\delta+(2\log a+2)\delta^2+\sum_{k=1}^\infty \frac{4(-1)^{k-1}}{(k+2)(k+1)k a^k}\delta^{k+2}=0.\]
\[\delta=\frac{1}{4\log a}-\frac{2\log a+2}{4a\log a}\delta^2-\sum_{k=1}^\infty \frac{(-1)^{k-1}}{(k+2)(k+1)k a^{k+1}\log a}\delta^{k+2}.\]
So that $\delta=1/4\log a+\Orden(|a|^{-1})$. 

Taking $R$ large enough, we then have $|b-a|<1$ with $a\in\Omega$ and $v(b)=2n$. It follows that $\Im v(b)=0$. The Lemma \ref{L:InOmega} implies that $b\in\Omega$. Proposition \ref{givesn} proof that in this case $n$ is a natural number and $b=\eta''_n$. 

If there were another point $a'\in\Omega$ with $u(a')=2n$, and we take $R$ large enough, then $|a-a'|<|a-\eta''_n|+|\eta''_n-a'|<1$, and Rouche's Theorem will give two solutions of $v(\eta)=2n$, contrary to Proposition \ref{givesn}. 
\end{proof}

\begin{proposition}
There is a natural number $N$ such that for each $n>N$ there is a zero $\rho_{-n}$ of $\Rzeta(s)$ in the fourth quadrant such that $|\rho_{-n}-2\pi i\overline\eta^2_n|<(2\pi|\eta_n|)^{-2}$, where $g(\eta_n)=1$. 
For each zero $\rho$ of $\Rzeta(s)$ in the fourth quadrant with $|\rho|>R$ there is some $n$ with $\rho=\rho_{-n}$. 
\end{proposition}
\begin{proof}
We apply Theorem \ref{X-rayTheorem} to the functions $f$, $g\colon\Omega\to\C$
\[f(\eta)=\frac{\chi(-2\pi i\eta^2 )\Rzeta(1+2\pi i\eta^2)}{\zeta(-2\pi i \eta^2)},\qquad 
g(\eta)=\exp\bigl(2\pi i \eta^2\log\eta-\pi i \eta^2+\pi i\eta-\tfrac12\log2-\tfrac{\pi i}{8}\bigr).\]
The conditions of Theorem \ref{X-rayTheorem} are satisfied. The region $\Omega$ is the one defined in Proposition \ref{P:defregion}. It is clear that $\Omega_r$ is connected for $r>0$.
We have seen (Proposition \ref{facts}) that $f(\eta)=g(\eta)(1+\Orden(|\eta|^{-1}))$ and we have 
\[\frac{g'(\eta)}{g(\eta)}=4\pi i\eta\log\eta+\pi i,\]
so that condition (b) in Theorem \ref{X-rayTheorem} is satisfied.

To prove condition (c) we notice that  if $f(\eta)=1$, then $\rho=2\pi i \overline{\eta}^2$ is a zero of $\Rzeta(s)$ in the fourth quadrant and then $\eta$ can be put in the form \eqref{rhoform}. Comparing with the definition of $\Omega$ in Proposition \ref{P:defregion} we see that for $R$ large enough the ball with center $\eta$ and radius $1$ is contained in $\Omega$. If $g(\eta')=e^{\pi i u(\eta')}=1$ with $|\eta'|$ large, then  Proposition \ref{nearaproxzeros} gives us a $\eta''$ with $|\eta'-\eta''|<1$ and such that $e^{\pi i v(\eta'')}=1$. Then Proposition \ref{givesn} gives a natural number $n$ such that $\eta''=\eta''_n$. If $R$ is large enough, we will have $n>N$ and $\eta''_n$ satisfies Proposition \ref{propeta''}. So $\overline D(\eta''_n,2)\subset\Omega$. This shows that the circle $\overline D(\eta',1)\subset \Omega$. Theorem \ref{X-rayTheorem} applies.

in what follows, we apply several propositions. In each one we can find a constant $R$ satisfying some properties. Here we take $R$ greater than all these different constants.

By Proposition \ref{Rzetaf}, there is $R$ such that if $\rho$ is a zero of $\Rzeta(s)$ in the fourth quadrant with $|\rho|>R$, then $\rho=2\pi i\overline\eta^2$ where $\eta\in\Omega$ satisfies $f(\eta)=1$. By Theorem \ref{X-rayTheorem} if $|\eta|>R$ there is some $\eta'$ with $g(\eta')=1$ with $|\eta-\eta'|<|\eta'|^{-2}$. Since $e^{\pi i u(\eta')}=1$ there is an integer $n$ with $u(\eta')=2n$. By Proposition \ref{nearaproxzeros} there is a solution of $v(\eta'')=2n$ such that $|\eta'-\eta''|<1/4\log|\eta'|+\Orden(|\eta'|^{-1})$and  $\eta''\in\Omega$. Then the Proposition \ref{givesn} shows that $n$ is a natural number and $\eta''=\eta''_n$.   Therefore, $\eta=\eta_n$ the unique solution in $\Omega$ of $u(\eta)=2n$. This shows that each zero $\rho$ is one of the $\rho_{-n}$. (We denote by $\rho_n$ for $n=1$, $2$, \dots\ the zeros of $\Rzeta(s)$ with $\Im(\rho_n)>0$. Therefore, the zeros in the fourth quadrant are denoted by $\rho_{-1}\approx10.64-i\;0.95$, $\rho_{-2}\approx16.43-i4.79$, \dots. If we start at the index $0$ we have to change the notation in this paper, so there is no zero $\rho_0$). 
\end{proof}

\section{Algorithm to compute the zeros}\label{S:6}

We have seen that $\rho$ is a zero of $\Rzeta(s)$ in the fourth quadrant if and only if $\rho=2\pi i\overline\eta^2$ with $\eta\in\Omega$ such that $f(\eta)=1$. By \eqref{Siegeleq} these are approximated by the zeros of 
\[S(\eta):=1+\eta^{2\pi i\eta^2}e^{-\pi i \eta^2}\frac{\sqrt{2} e^{\frac{3\pi i}{8}}\sin(\pi\eta)}{2\cos(2\pi\eta)}.\]
This is very nice, as can be seen in the X-ray of $S(\eta)$ in Figure \ref{second} below, where we have added a little dot at every point 
\[\eta_n=\sqrt{\frac{i\overline{\rho}_{-n}}{2\pi}},\]
where $\rho_{-n}$ are the zeros of $\Rzeta(s)$ in the fourth quadrant ordered by the absolute value. This means that in practice our theorems are valid for small $R$ and small $N$.

\begin{algorithm}
	\caption{Computing $\rho_{-n}$} 
	\begin{algorithmic}[1]
	 \State Compute $a_n=\frac{1}{8}-2n-\frac{\log2}{2\pi}i$
	 \State  Compute $\eta_n''$ from $\eta_n''^2 = a_n/W_1(a_n/e)$ with $\pi/4<\arg(\eta'')<\pi/2$
	  \State Use $\eta''_n$ as an initial point to compute the solution $\eta'_n$ of 
	  $\eta^2\log\eta^2-\eta^2+\eta=a_n$
      \State Compute $\tilde\rho_n$,  solution of $\Rzeta(s)-\chi(s)\zeta(1-s)=0$ starting from the approximate solution $\tilde\rho' = 1+2\pi i (\eta'_n)^2$
      \State Extract the zero  $\rho_{-n}=1-\overline{\tilde\rho_n}$
	\end{algorithmic} 
\end{algorithm}
Let $a_n=\frac{1}{8}-2n-\frac{\log2}{2\pi}i$, where $n$ is a natural number. 
Let $W:=W_1(a_n/e)$. By \eqref{eqWeta} we have $\log(\eta^2/e)=W$. By the definition of the Lambert function $We^W=a_n/e$, or $e^W=\frac{a_n}{eW}$. Therefore, ${\eta''_n}^2=a_n/W$.

The next table gives an idea of the approximation obtained for small values of $n$,
\[\begin{array}{cccc}
n & 1 & 100 & 1000 \\ \hline
\noalign{\smallskip}
\eta''_n  & 1.022367+i\; 0.644959  & 2.447556+i\; 6.789936 &   4.658712+i\; 18.391214\\
\eta'_n & 0.875036+i\; 0.967592 & 2.355945+i\; 6.846082 &   4.587993+i\; 18.422872\\
\eta_n   & 0.880367+i\; 0.962502 & 2.355977+i\; 6.846003 &   4.587996+i\; 18.422864\\
\end{array}\]

A table with the first 2122 zeros in the fourth quadrant is included in the \texttt{.tex} file of this document. They are computed with the algorithm presented in this section. The real and imaginary parts of each zero are given with 25 correct decimal digits.

\begin{figure}[H]
\begin{center}
\includegraphics[width=0.9\hsize]{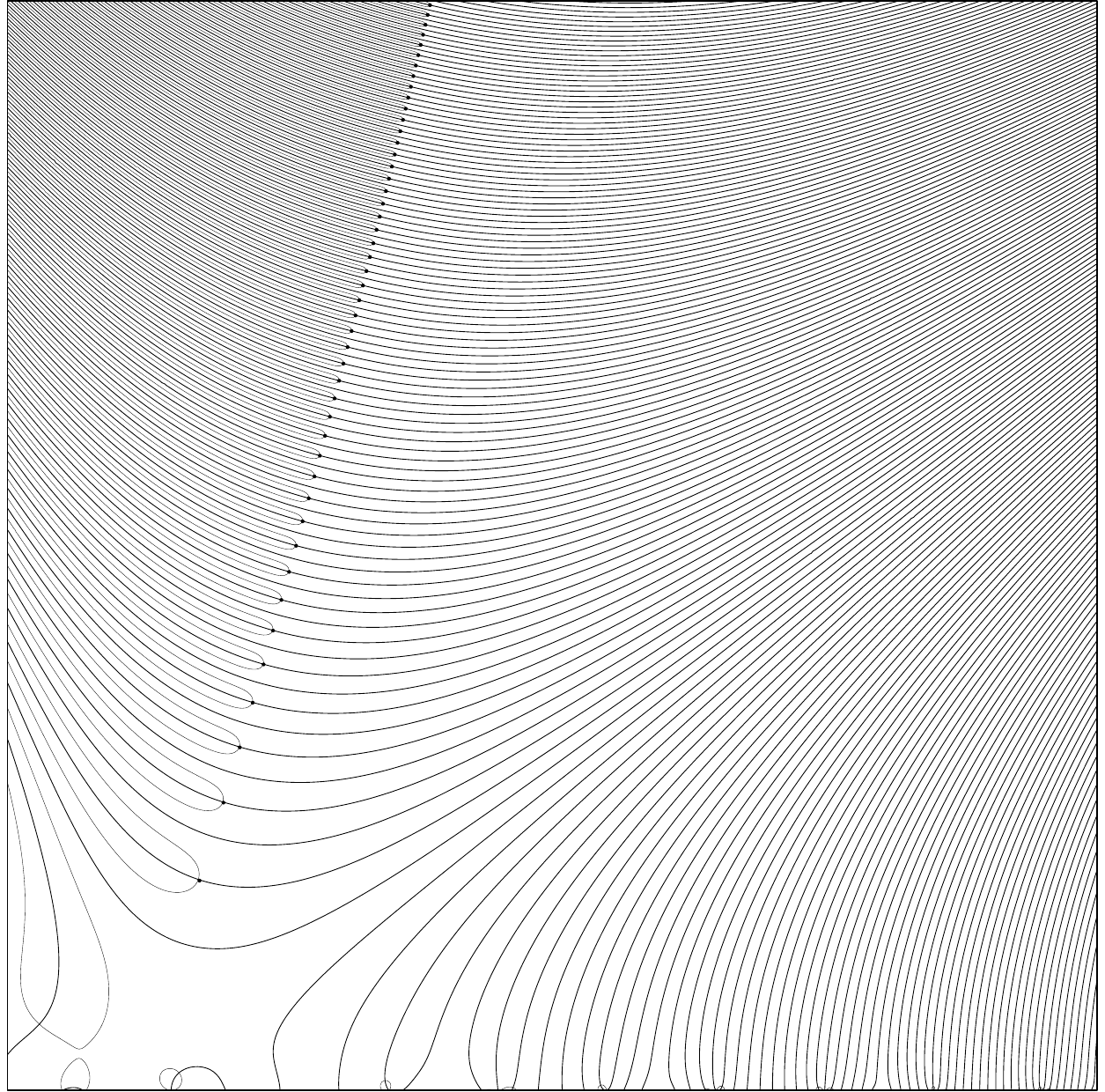}
\caption{X-ray of $S(\eta)$ on the rectangle $(0,5)^2$ and dots at the $\eta_n$.}
\label{second}
\end{center}
\end{figure}

\section{Asymptotic of zeros}\label{S:7}

\begin{proposition}
The zero $\rho_{-n}$ of $\Rzeta(s)$ in the fourth quadrant has an asymptotic expansion whose first terms are given in \eqref{expansion}. In particular $\rho_{-n}=\beta_{-n}+\gamma_{-n}$ with 
\begin{equation}
\beta_{-n}\sim \frac{4\pi^2n}{\log^2n},\qquad \gamma_{-n}\sim-\frac{4\pi n}{\log n}.
\end{equation}
\end{proposition}
\begin{proof}
For $n$ large, the number $\eta''_n$ satisfies the equation \eqref{eqWeta}. 
By \cite{C} we get an asymptotic expansion for  $\eta''_n$ as in \eqref{firstasymptotic}. After substituting this expansion into \eqref{E:toasymp}, and eliminating terms of smaller order than any negative power of $(\log n)$ we obtain the relation 
\begin{equation}\label{eq}
-\frac{2}{x}\Bigl(1+\sum_{k=1}^\infty\frac{U_k(y)}{x^k}\Bigr)^2\Bigl\{A+x-y +2\log\Bigl(1+\sum_{k=1}^\infty\frac{U_k(y)}{x^k}\Bigr)-1\Bigr\}+2=0,
\end{equation}
where $A=\pi i +\log 2$, $x=\log n$ and $y=\log\log n$. Expanding \eqref{eq} to powers of $x^{-1}$ and equating the coefficients to $0$ we get the values of the polynomials $U_k(y)$.

Therefore, we get 
\[\eta''_n=i\sqrt{\frac{2n}{\log n}}\Bigl(1+\sum_{k=1}^K\frac{U_k(\log\log n)}{\log^k n}+\Orden(\log^{-K-1}n)\Bigr).\]
Since $|\eta''_n-\eta'_n|\le 1/4\log|\eta_n'|+\Orden(|\eta_n'|^{-1})$, the same approximation is true for $\eta'_n$ (note that the error in the above equation for $\eta''_n$ is of order 
$\Orden(n^{1/2}\log^{-K-3/2}n)$). In the same way, this same expansion is true for $\eta_n$. Since $\rho_{-n}=2\pi i\overline\eta_n^2$, we get an expansion for $\rho_{-n}$ of the form
\[\rho_{-n}=-2\pi i \frac{2n}{\log n}\Bigl(1+\sum_{k=1}^K\frac{\overline{U_k(\log\log n)}}{\log^k n}+\Orden(\log^{-K-1}n)\Bigr)^2.\]
Expanding the square and eliminating terms greater than the error, we get the asymptotic for $\rho_{-n}$.

The first  terms of this expansion are (with $\ell_n=\log\log n$)
\begin{equation}\label{expansion}
\begin{aligned}
\rho_{-n}&=\frac{4\pi^2 n}{\log^2 n}\Bigl(1+\frac{2\ell_n-2\log2-1}{\log n}+\frac{3\ell_n^2-(5+6\log 2)\ell_n+3\log^2 2+5\log 2+1-\pi^2}{\log^2n}\Bigr)\\
&-i\frac{4\pi n}{\log n}\Bigl(1+\frac{\ell_n-\log 2}{\log n}+\frac{\ell_n^2-(2\log 2+1)\ell_n+\log^22+\log2-\pi^2}{\log^2n}\Bigr)\\
&+\Orden(n(\log\log n)^3/(\log n)^4).
\end{aligned}
\end{equation}
\end{proof}

\end{document}